\documentclass[12p]{amsart}
\usepackage{color}
\usepackage{graphicx}

\usepackage[colorlinks=true]{hyperref}
\hypersetup{urlcolor=blue, citecolor=blue, draft=false}
\newtheorem{definition}{Definition}[section]
\newtheorem{Theorem}{Theorem}[section]

\theoremstyle{remark}

\theoremstyle{remark}

\begin{document}

\title[ON BUNGEE SETS OF COMPOSITE TRANSCENDENTAL ENTIRE FUNCTIONS]{ON BUNGEE SETS OF COMPOSITE TRANSCENDENTAL ENTIRE FUNCTIONS }

\author[Anand P. Singh ]{ Anand P. Singh}
\address{Department of Mathematics, Central University of Rajasthan,  Bandarsindri, Kishangarh-305817, Rajasthan, India}
\email{ singhanandp@rediffmail.com }

\date{}

\bigskip
\begin{abstract} Let $f$ be  a transcendental entire function. For $n \in \mathbb{N},$ let $ f^{n}$ denote the $n^{th}$ iterate of $f$. Let 
$ I(f) = \{z \in \mathbb{C} : f^n \rightarrow \infty $  as $  n \rightarrow \infty \} $ and  
$ K(f) = \{z: \textrm{ there exists } R > 0 \textrm{ such that } | f^n(z) | \leq R  \textrm{ for  }   n \geq 0 \}. $ Then the set $ \mathbb{C}\ \setminus (I(f) \cup K(f)) $  denoted by $ BU(f) $ is called Bungee set of $f$.  In this paper we study some properties of Bungee sets of composite transcendental entire functions and also of Bungee sets of permutable transcendental entire functions.

\end{abstract}

\keywords{Transcendental Entire Function, Escaping set, Julia set, Bungee set \\
 Sponsored by SERB grant No. EMR/2016/001645}

{\Large
\bigskip
\maketitle

\section{Introduction}
\setcounter{equation}{0}

Let $f$ be a rational or a transcendental entire function. For $n \in \mathbb{N} $, let $ f^{n}$ denote the $n^{th}$ iterate of $f$. The set \\
$ F(f) = \{ z  : (f^n) $ is normal in some neighbourhood of $ z \}$  is called Fatou set of $f$ and its complement denoted by $J(f)$ is called Julia set of $f$. For properties of these two sets one can refer for instance \cite{Beardon}, \cite{Bergweiler}, \cite{Morosawa}. Here we observe that  $\hat{\mathbb{C} }= \mathbb{C} \cup \{ \infty \}$ has been partitioned in two different classes of sets, viz. the Fatou set and the Julia set with  $ \infty  \in J(f) $ when $f$ is transcendental entire. There are other ways also that $\hat{\mathbb{C} }$ can be partitioned, which is also linked with the Fatou and Julia sets. 

For a transcendental entire function, Eremenko \cite{Eremenko} considered  the Escaping set 
$ I(f) = \{z \in \mathbb{C} : f^n(z) \rightarrow \infty \textrm{ as } n \rightarrow \infty \}  $
and showed that $ I(f) \cap J(f) \not= \phi $, $ \partial I(f) = J(f) $ and all components of $ \overline{I(f)}$ are unbounded. He further conjectured that all the components of $I(f)$ should be unbounded. This conjecture though not yet solved has given rise to a rich development in the field. The points which tend to infinity with different ``speeds"  such as fast escaping points, slow escaping points, relatively fast escaping points etc. have been studied by various authors,  (see for instance \cite{Rippon}, \cite{Rippon 1}, \cite{ Singh}, \cite{ Wang},  \cite{Zheng}).  In contrast, there are points whose iterates under $f$ remain bounded. The set of such points is denoted by  $ K(f) $ and is defined by 

    $ K(f) = \{z: \textrm{ there exists } R > 0 \textrm{ such that } | f^n(z) | \leq R  \textrm{ for  }   n \geq 0 \}. $
    
 For non linear polynomial $P,$ $ K(P) $ is known as filled Julia set, which has been extensively studied. For transcendental entire function, $K(f) $ has been studied for instance  by Bergweiler \cite{ Bergweiler 1} and Osborne \cite{Osborne}. There are points whose iterate are neither bounded nor tending to $\infty$. Such sets have been considered by Osborne and Sixsmith. They called it Bungee set. More specifically Osborne and Sixsmith \cite{Osborne 1}, defined the Bungee set $ BU(f) $ as $ BU(f) = \mathbb{C} \setminus  ( I(f) \cup K(f) )$. Note that as each of   $I(f)$ and $K(f)$ are  completely invariant, it follows that $ BU(f)$ is also completely invariant, and the sets $I(f)$, $K(f)$, and $ BU(f)$ partition the complex  plane $ \mathbb{C} $  into three classes. Rather than studying the Bungee set as complement of $ I(f) \cup K(f) $  we give here an alternate  definition for Bungee set, which is very easy to use and prove some results using this definition.
 
\begin{definition} Let $f$ be rational or transcendental entire function. We define the Bungee set of $f$ denoted by $BU(f)$ by: 
 $ BU(f) = \{ z :  \textrm{ there exist atleast two subsequences } \{f^{n_k}\},  \{f^{m_k}\}  \\ \textrm{ with } {n_k} \rightarrow \infty \textrm{ and  } {m_k} \rightarrow \infty  \textrm{  as } k   \rightarrow \infty  \textrm{ and a constant  } R > 0 \textrm{    such }  \textrm{that } |f^{n_k}(z) | \leq R  \textrm{ for } k =1, 2, \dots,  \textrm{ and }$  $   f^{m_k}(z)  \rightarrow \infty  \textrm{ as } m_k \rightarrow \infty  \}. $ 
 \end{definition}
 
 Note that the  subsequences in the above definition  cannot have a ``pattern" in the sense that for the point $z_o \in BU(f) $ there do not exist sequences  $ {n_k}$ and  $ {m_k} $    such that $ n_k = k n_o, m_k = k m_o, $ where  $ n_o, m_o \geq 1, k = 1, 2, \dots $  and such that  $ |f^{n_k}(z_o) | \leq R  \textrm{ for } k =1, 2, \dots,  \textrm{ and }   f^{m_k}(z_o )  \rightarrow \infty  \textrm{ as } k \rightarrow \infty $, for  if they do exist, then choose $k$ sufficiently large say $l$ such that $  |f^{m_k}(z_o) | > R$  for all $ k \geq l $, and so $  |f^{k m_o}(z_o) | > R$ for all $ k \geq l $, and in particular   $  |f^{k m_o n_o}(z_o ) | > R $ for   $ k \geq l $, contradicting $  |f^{t n_o}(z_o ) | =  |f^{n_t}(z_o ) | \leq R$ for all $t$.  Also, note that  for $z_o \in BU(f)$ if $ ( n_k, m_k, R, l_o )$ denote the sequences $ \{n_k\}$ and  $ \{m_k\} $    such that  $ |f^{n_k}(z_o) | \leq R  \textrm{ for } k =1, 2, \dots,  \textrm{ and }  |  f^{m_k}(z_o ) | > R  $, for all $ k \geq l_o $, then if 
 $ n_k = p  n_{k-1} $ and  $ m_k  = q  m_{k-1},  $ then $ p \not= (  \frac{m_o}{n_o} q^l )^{\frac{1}{k}}$ for all $ l \geq l_o. $ For if  for some $ l \geq l_o $ $ p = (  \frac{m_o}{n_o} q^l )^{\frac{1}{k}}$, then $ n_k = p^k  n_o = q^l m_o = m_l, $ and so  $ |f^{n_k}(z_o) | = |f^{m_l}(z_o) | $ which is not possible as left equality is bounded by $R$ and the right equality is greater than $R$.
  \\

    For a non-linear polynomial $P, BU(P) = \phi $, and for a transcendental entire function $f$,  $BU(f) \cap J(f) \not= \phi $ \cite{Baker 1}, and also there are examples where  $BU(f) \cap F(f) \not= \phi $ (see for instance \cite{Bergweiler 1}, \cite{Bishop}). For rational functions, the Bungee set may coincide with Fatou set, for instance if $R(z) = \frac{1}{z^2}$ then $BU(R)
= (|z|<1)\cup (|z|>1) $ and $J(R) = \partial (BU(R))$. It is interesting to note that if $R$ and $S$ are permutable rational functions  ( i.e., $ R \circ S = S \circ R$), then $J(R)= J(S)$, where as this is still an open question for transcendental entire function, though quite some progress has been made (see for instance \cite{Baker}, \cite{Singh}). However the corresponding result on Bungee set   is not true for rational functions. For instance if $R(z) = \frac{1}{z^2} $ and $S(z) = z^2 $ then clearly $ R(S(z)) = S(R(z)) $, and $ BU(R)= (|z|<1) \cup (|z|>1)  \textrm{ whereas } BU(S) = \phi $. It looks, the same type of statement would be true if $f$ is transcendental entire, though we do not have an example for it. Thus the Bungee set may behave quite differently  from Julia sets for rational functions and transcendental entire functions. It is also well known that for any $n$, $J(f^n) = J(f),  F(f^n) = F(f).$ However this is not true for Bungee sets. We only have $ BU(R^n) \subset BU(R)$, since for instance, if $R(z) = \frac{1}{z^2}$, then  $ BU(R^2) = \phi , \textrm{ whereas }   BU(R)= (|z|<1) \cup (|z|>1)$. We are not sure whether  the equality  holds for transcendental entire function. But we can  atleast show $ BU(f^n) \subset BU(f)$. For this observe that if $z_o \in BU(f^n)$, then there exist $\{n_k\}$ and $\{m_k\} $ tending to $ \infty$ such that $|(f^n)^{n_k}(z_o)| \leq R $ for some positive $R$, $ k = 1, 2, \dots$ and $(f^n)^{m_k}(z_o) \rightarrow \infty $
as  $ m_k \rightarrow \infty $. Let $ t_k = n \cdot n_k $ and $p_k = n \cdot m_k$, then  $|(f)^{t_k}(z_o)| \leq R ,  k = 1, 2, \dots$ and $f^{p_k}(z_o) \rightarrow \infty $ as  $ p_k \rightarrow \infty $, and so $z_o \in BU(f)$.\\
\newpage

\section{Main results and their proofs.}
\setcounter{equation}{0}

  In this paper we concentrate on Bungee sets of transcendental entire functions, and more specifically with  composition of transcendental entire functions.\\

If $U$ is a component of  $F(f)$, then by complete invariance, $f(U)$ lies in some component of $F(f)$. If $U_n \cap U_m = \phi $ for $ n \not= m$ where $U_n$ denotes the component of $F(f)$ which contains $f^n(U)$, then $U$ is called a wandering domain, else it is either a periodic or a pre-periodic domain. A complete classification for periodic domain is well known ( see for instance \cite{Morosawa} ). Also it is well known \cite{Sullivan} that rational functions have no wandering domain,  though the same is not true for transcendental entire functions. Infact transcendental entire functions may  have wandering domains. Also there are transcendental entire functions, such as the functions in Speiser class ( i.e. entire functions whose singularities are finite in number)
which do not have wandering domain. There are other classes of transcendental entire functions which do not have wandering domain (see for instance  \cite{Bergweiler}, \cite{Morosawa}).

If $P$ is a bounded periodic component of $F(f)$ then obviously every point $z_0 \in P$ does not belong to $BU(f)$. But what can be said about points on $\partial P$ ? We shall show that even $ \partial P$ does not contain points of $BU(f)$. Thus we shall prove

\begin{Theorem}\label{Theorem 1}
  Let $f$ be transcendental entire function. Let   $P$ be a bounded  periodic component of $ F(f)$. Then $ \partial P \cap BU(f) =  \phi $.
 \end{Theorem}

\noindent \textit{Proof of Theorem} \ref{Theorem 1}. Note that as $P$ is a bounded periodic Fatou component, there exists a constant $M$ such that $ |f^n(P)| \leq M $ for all $ n \in \mathbb{N}$. 
Suppose $ \xi_o \in  \partial P \cap BU(f).  $
Then  there exists a sequence $\{ m_k\}$ such that $ f^{m_k}(\xi_o) \to \infty $ as $ m_k \rightarrow \infty. $ 
Hence for sufficiently large $m_k$ say $ M_k$, $ | f^{M_k}(\xi_o)\mid  >  1+M.$ Also $ f^{M_k}$ is continuous at 
$\xi_o$  and hence there exists $ \delta > 0$ such that $ | z - \xi_o |  <  \delta $ implies 
$ | f^{M_k}(z) - f^{M_k}(\xi_o) |  < 1$. In particular for $ t  \in  ( | z - \xi_o |  < \delta )\cap P $ we have $ | f^{M_k}(t) - f^{M_k}(\xi_o) |  < 1$, and on the other hand 
$ | f^{M_k}(t) - f^{M_k}(\xi_o) | \geq  | f^{M_k}(\xi_o) |- |f^{M_k}(t) | >1. $ This contradiction proves the theorem.\\

{\bf Question}: Is $ \partial P \cap BU(f) = \phi $ if $P$ is unbounded Fatou component ? \\

We next show that  atleast for  a transcendental entire function with no wandering domain, the bungee set cannot be a closed set.

\begin{Theorem}\label{Theorem 2}
  Let $f$ be transcendental entire function without wandering domain. Then  $  BU(f) $ cannot be a closed subset  of $ \mathbb{C} $. 
 \end{Theorem}
 
 For the proof of the theorem we shall need the following theorem of 
 Osborne and Sixsmith  \cite{Osborne}
 
\begin{Theorem}\label{Theorem 3}
Let $f$ be a transcendental entire function.\\
(a) If $U$ is a Fatou component of $f$ and $U \cap BU(f) \not= \phi,$ then $ U \subset BU(f) $ and $U$ is a wandering domain.\\
(b) $J(f) = \partial BU(f).$\\

\end{Theorem}

\noindent \textit{Proof of Theorem} \ref{Theorem 2}. Suppose $  BU(f) $ is a closed subset  of $ \mathbb{C}. $  As $  BU(f) \not= \phi $ it contains infinitely many points. Also clearly $ BU(f)  $ is completely invariant set. Also as  $J(f)$ is the smallest closed completely  invariant set having atleast three points it follows that $ J(f) \subset BU(f)$. If $ J(f) \not=  BU(f)$, then there exists $ z_o \in BU(f) \cap  F(f) $, and consequently there exists $ U \subset F(f) $ which intersects $BU(f)$. By Theorem \ref{Theorem 3}, $U$ must be a wandering domain. This contradiction proves $ J(f) =  BU(f).$ But then this result is also not possible as repelling periodc points lie in $J(f) $ and periodic points obviously do not lie in $BU(f)$. This completes the proof of the  Theorem.\\

Our next two theorems deal with pre-images of bungee sets corresponding to permutable transcendental entire functions. 

\begin{Theorem}\label{Theorem 4}
 Let $f$ and $g$ be permutable transcendental entire functions. Let $U \subset BU(f).$ If 
$g^{-1} (U) \not= \phi $, then  $g^{-1} (U) \cap ( I(f) \cup BU(f) ) \not= \phi $\\
 \end{Theorem}

\noindent \textit{Proof of Theorem} \ref{Theorem 4}. Suppose $ g^{-1} (U) \cap ( I(f) \cup BU(f) ) = \phi $. Let $z_o \in g^{-1}(U) $. Then there exists a constant $A$ such that $  |f^n(z_o))| \leq A $ for all  $n \in \mathbb{N} $, and consequently $ | g(f^n(z_o))| \leq M(A, g) $ for all $ n \in \mathbb{N}$. Also there exists $\xi_o \in U$ such that $ g(z_o) = \xi_o$, and as 
$ U \subset BU(f)$ there exists a sequence $n_k$ such that $f^{n_k}(\xi_o) \rightarrow \infty $ as $ n_k \rightarrow \infty$. So choose $n_k$ sufficiently large so that $ |f^{n_k}(\xi_o) | > M(A,g) +1$. Then for such $n_k$, $ |f^{n_k}(g(z_o) )| > M(A,g) +1$ and where as $ |f^{n_k}(g(z_o) )| = |g(f^{n_k}(z_o) ) | \leq  M(A,g) .$
 This contradiction proves the theorem.\\
 
 \begin{Theorem}\label{Theorem 5}
 Let $f$ and $g$ be permutable transcendental entire functions. Further let there exist a non linear polynomial $P$ such that $ P\circ f = f \circ g $.  If  $U \subset BU(f), $ 
 then  $g^{-1} (U)   \subset BU(f).$
 
  \end{Theorem}
 
\noindent \textit{Proof of Theorem} \ref{Theorem 5}. Following as in Theorem   \ref{Theorem 4},  if $z_o \in g^{-1}(U) $ and if $  |f^n(z_o))| \leq A $ for all  $n \in \mathbb{N} $, then we have obtained a contradiction. We next show that $ f^n(z_o) \not \rightarrow \infty $ as $ n \rightarrow \infty$, and so $z_o \in BU(f)$. So, let $ f^n(z_o) \rightarrow \infty $ as $ n \rightarrow \infty$. Then there exists $\xi_o \in U $ and a constant $K$  and a sequence $n_k \rightarrow \infty $ such that $  | f^{n_k}(\xi_o)| < K $ for all $ k=1,2,
\dots $. Let $ p_k = n_{k+1}-n_k $ for all $k$. As $P(z)$ is a non linear polynomial, there exists $R$ sufficiently large such that $|P(z)| >  |z| $ for all $z \in (|z| >R) $.As $ f^n(z_o) \rightarrow \infty $ as $ n \rightarrow \infty$, we can choose $n_k $ sufficiently  large so that\\
 $  |f^{n_k}(z_o) | > Max \{R, K+1\} , k=1,2,\dots .$  Let $ f^{n_k}(z_o) = \eta $. Then $ f^{n_k}(g(z_o)) = f^{n_k}(\xi_o)= t $ say where $|t|< K$. And so by permutability, $ t = g(f^{n_k}(z_o))= g(\eta) $, and so 
 $ | f^{p_k}(g(\eta ))| = | f^{p_k}(g(f^{n_k}(z_o)))| = | f^{p_k}(f^{n_k}g(z_o))|= | f^{n_{k+1}}(g(z_o))| = | f^{n_{k+1}}(\xi_o))| < K $.\\
 On the other hand, $| f^{p_k}(g(\eta))| =  | f^{n_{k+1}}(g(z_o))|  = | P(f^{n_{k+1}}(z_o))| \geq | (f^{n_{k+1}}(z_o)| > K + 1$, and so  $ z_o \in BU(f)$.  Thus $g^{-1}(U) \subset BU(f)$.\\

  If $f$ is a transcendental entire function and  $z_o \in BU(f)$, then there exists $R > 0 $ and sequences 
$\{n_k\}$ and $\{m_k\} $ tending to $ \infty$ such that $ | f^{n_k}(z_o)| \leq R $ for some positive $R$, $ k = 1, 2, \dots$ and $ f^{m_k}(z_o) \rightarrow \infty $
as  $ m_k \rightarrow \infty $. We shall denote such domain $  ( |z| \leq R ) $ by $D_R(z_o)$. (Note that such $D_R(z_o)$ is not unique).\\

\begin{Theorem}\label{Theorem 6}
 Let $f$ and $g$ be permutable transcendental entire functions. Let $P$ be a polynomial of degree $\geq 2$ and let $h$ be  a transcendental entire function such that $ P\circ f = h \circ g$. Let $z_0 \in BU(f)$. Then $ g(D_R(z_o))$ cannot be a periodic domain of $f$.

 \end{Theorem}

\noindent \textit{Proof of Theorem} \ref{Theorem 6}. Let $z_o \in BU(f)$. Then there exists $R > 0 $ and sequences 
$\{n_k\}$ and $\{m_k\} $ tending to $ \infty$ such that $| f^{n_k}(z_o)| \leq R $ for positive $R$, $ k = 1, 2, \dots$ and $ f^{m_k}(z_o) \rightarrow \infty $
as  $ m_k \rightarrow \infty.$\\
First suppose $ |f^n(g(D_R(z_o)))| \leq A $ for some constant $A$ and all $ n \in \mathbb{N}$. Then 
$  |h(f^n(g(D_R(z_o))))| \leq M(A,h) < B $ say, for all $n \in \mathbb{N}$ where $ M(A,h) = \max _{|z|=A}|h(z)|$.
As $P$ is a polynomial  of degree $ \geq 2 $, we can choose $S$ sufficiently large so that $ |P(z)|>|z|$ for all $|z|>S$.\\
Now choose $M_k$ sufficiently large so that $ | f^{m_k}(z_o)| > \max \{S, B \}  $ for all $ m_k \geq M_k.$\\
Thus  for  $ m_k \geq M_k, $ $ | P( f^{m_k}(z_o))| > | f^{m_k}(z_o)| > B.$\\
On the other hand $ | P( f^{m_k}(z_o))| =  | P( f^{m_k- n_1}(f^{n_1}(z_o)))| = | P( f^{m_k-n_1}(t))| = | h( f^{m_k-n_1-1}g(t))| \leq M(A,h) < B $ where $t= f^{n_1}(z_o)$. This contradiction shows that $ f^n(g(D_R(z_o))) $ cannot be bounded for all $n \in \mathbb{N}$ and hence $g(D_R(z_o))$ cannot be a periodic Fatou component of $f$ unless it is a periodic Baker domain, which also is not possible as  $g(D_R(z_o))$ is bounded.\\

If $f$ and $g$ are permutable transcendental entire functions and if $z_0 \in BU(f)$ as well as in $F(f)$ then we have the following result.

\begin{Theorem}\label{Theorem 7}
 Let $f$ and $g$ be permutable transcendental entire functions. Let $h$ be a trascendental entire function and $P$ be a non-linear polynomial such that $ P \circ f = h \circ g $. If $ z_0 \in F(f) \cap BU(f)$ then $ g^{-1}(z_0) \in F(f) \cap BU(f)$.

 \end{Theorem}
 
\noindent \textit{Proof of Theorem} \ref{Theorem 7}. Since $ z_0 \in F(f) $ and $F(f)$ is open, there exists a neighbourhood and consequently a component $U$ of Fatou set of $f$ such that $U \cap BU(f) \neq \phi$, and so by Theorem \ref{Theorem 3}, $U \subset BU(f)$ and $U$ is a wandering domain and $\partial BU(f) = J(f)$. If $ g^{-1}(z_0) \in J(f)$, then by complete invariance $ g( g^{-1}(z_0)) \in J(f)$. Thus $z_0 \in J(f)$, contradicting $z_0 \in F(f).$

We now show that $ g^{-1}(z_0) \in BU(f)$. For suppose $ |f^n(g^{-1}(z_0)| \leq A  $ for some constant $A$ and for all $n \in \mathbb{N}$, then $  |g(f^n(g^{-1}(z_0)) | \leq M(A, g) $, and so as $f$ and $g$ are permutable,    $  |(f^n(z_0)) | \leq M(A, g) $ for all $n\in \mathbb{N}$ contradicting $ z_0 \in  BU(f)$. Next suppose  $ | f^n(g^{-1}(z_0) | \rightarrow \infty,$ as $ n \rightarrow \infty. $ Now, as $ z_0 \in  BU(f)$, there exist  a subsequence 
 $\{n_k \} \textrm { tending to }\infty $  and a constant  $ c > 0 $,  such that $ | f^{n_k}(z_0)| \leq c, k=1,2,\dots. $   
Let $ D = M(c, h)$. Since $ P $ is a nonlinear polynomial, we  can choose $ r > D+1$ sufficiently large so that  for all $ z \in (|z|>r),   |P(z)| > |z|$. Since $ | f^n(g^{-1}(z_0) | \rightarrow \infty $  as $ n \rightarrow \infty $,  we select $n_t$ from the subsequence so that

   $ | f^{n_t}(g^{-1}(z_0) | > r $ and   $ | f^{n_t+1}(g^{-1}(z_0) | > r $.
   
  Now$ (f^{n_t}\circ g)(g^{-1}(z_0) = f^{n_t}(z_0) = \xi_0 $ say, where $ |\xi_0 | \leq c $.\\
  And so\\
  $ \xi_0 = (f^{n_t}\circ g)(g^{-1}(z_0) = g(f^{n_t}(g^{-1}(z_0)) = g(\eta ) $  say, where $|\eta| = |f^{n_t}g^{-1}(z_0)| > r$.
  Now $ | (h \circ g) (\eta )| = |h(\xi_0)|  \leq M(c, h) = D$, and also  as $ |f^{n_t+1}(g^{-1}(z_0) | > r $, we have\\
  $|(h \circ g) (\eta )| = | (P\circ f)(\eta) | = | P( f^{n_t+1}(g^{-1}(z_0)) | > | f^{n_t+1}(g^{-1}(z_0)) | > r > D+1. 
   $  This contradiction proves the theorem.\\

\begin{Theorem}\label{Theorem 8}

 Let $f$ be a transcendental entire function. If $ z_0 \in BU(f \circ g)$ then $ g(z_0) \in BU(g \circ f)$.
 
 \end{Theorem}
 \noindent \textit{Proof of Theorem} \ref{Theorem 8}. Suppose $ g(z_0) \notin BU(g \circ f)$. Then  either (i) $ | (g \circ f )^n (g(z_0)) | \leq A $ for all 
$n \in \mathbb{N} $ and some $ A > 0$ or (ii) $ \lim_{n \rightarrow \infty }(g\circ f)^n (g(z_0)) = \infty $. \\
Now if (i) holds, then $ | g(f\circ g)^n(z_0) | \leq A $ for all  $n \in \mathbb{N} $, and so $ | f(g(f\circ g)^n(z_0)|\leq M(A, f) $ for all $ n \in \mathbb{N},$ where $ M(A, f) = \max_{|z|=A}|f(z)| $. Thus $ | (f\circ g)^{n+1}(z_0) | \leq M(A, f) $ for all $n\in \mathbb{N}$, contradicting $ z_0 \in  BU(f \circ g).$\\
Next, if (ii) holds, then since $z_0 \in BU(f \circ g)$, there exists a sequence ${n_k}$ such that $ |(f \circ g)^{n_k}(z_0) | \leq \beta,$ for some  $\beta > 0 $. And so $ | g(f\circ g)^{n_k}(z_0)| \leq M( \beta, g), k = 1, 2, \dots $. Thus $ | (g\circ f)^{n_k}(g(z_0))| \leq M( \beta , g)$, contradicting (ii). This proves the theorem.\\

     \end{document}